# RATE-OPTIMAL ESTIMATION FOR A GENERAL CLASS OF NONPARAMETRIC REGRESSION MODELS WITH UNKNOWN LINK FUNCTIONS

By Joel L. Horowitz[1] and Enno Mammen[2]

### Northwestern University and University of Mannheim

This paper discusses a nonparametric regression model that naturally generalizes neural network models. The model is based on a finite number of one-dimensional transformations and can be estimated with a one-dimensional rate of convergence. The model contains the generalized additive model with unknown link function as a special case. For this case, it is shown that the additive components and link function can be estimated with the optimal rate by a smoothing spline that is the solution of a penalized least squares criterion.

**1. Introduction.** This paper presents a general class of nonparametric regression models with unknown link functions. The models include neural network structures where link functions enter into the model on different levels. The inputs into the nodes of the net are modeled as sums of transformations of lower level inputs. Different approaches to modeling the transformations are allowed, including smooth nonparametric functions, shape-restricted nonparametric functions and parametric specifications. We show that rate optimal estimation in this class of models can be achieved by penalized least squares. The proof of the result relies on direct application of empirical process theory.

The approach described in this paper permits a unified treatment of a large class of models that includes some well-known examples. The proposed estimation method can be implemented in practice by using smoothing splines.

Received August 2005; revised January 2007.

[1]Supported in part by NSF Grant SES-03-52675 and the Alexander von Humboldt Foundation.

[2]Supported by Deutsche Forschungsgemeinschaft MA 1026/7-2.

The simplest form of our model is a generalized additive model with an unknown link function. That is,

$$Y = F[m_1(X^1) + \cdots + m_d(X^d)] + U, \tag{1}$$

where $X^1, \ldots, X^d$ are one-dimensional components of a $d$-dimensional covariate vector, $F$ and $m_1, \ldots, m_d$ are unknown functions and $U$ is an unobserved error variable satisfying $E[U|X] = 0$. We first discuss estimation of this model when all the unknown functions belong to the same smoothness class. We will show that these functions can be estimated with $L_2$-rate $n^{-k/(2k+1)}$ if they are $k$-times differentiable. Penalized least squares estimators with properly chosen penalty functions achieve this rate. The rate is optimal because it would be optimal if the link function were known. As a corollary, we will get the result that this rate carries over to models that assume more structure on $F$ and $m_1, \ldots, m_d$. Empirical process theory is our main tool for obtaining rate optimality. See van de Geer [40] for a comprehensive exposition of the use of empirical process theory in nonparametric estimation. Applying these techniques, it can be shown relatively directly that the function $(x_1, \ldots, x_d) \rightsquigarrow F[m_1(x_1) + \cdots + m_d(x_d)]$ can be estimated with rate $n^{-k/(2k+1)}$. The main difficulty is to show that this rate carries over to the estimation of the functions $F$ and $m_1, \ldots, m_d$. Clearly, identification of these functions requires normalizing restrictions.

If the link function, $F$, is known to be the identity function, then (1) is a nonparametric additive regression model. This model has been extensively studied. Stone [35, 36, 37] and Newey [30] have shown that optimal $L_2$-rates can be achieved by piecewise polynomial fits and regression splines. Breiman and Friedman [4] and Buja, Hastie and Tibshirani [5] discuss backfitting for additive models. Opsomer and Ruppert [34] and Opsomer [33] considered pointwise asymptotic distribution theory for backfitting. Mammen, Linton and Nielsen [22] introduced smooth backfitting estimates, a modification of backfitting that works more reliably in the case of many components and irregular design and that allows a complete asymptotic theory. Nielsen and Sperlich [31] and Mammen and Park [24, 25] discuss practical implementation of smooth backfitting. Tjøstheim and Auestad [38], Linton and Nielsen [21] and Fan, Härdle and Mammen [9] discuss marginal integration estimators. See Christopeit and Hoderlein [6] for a related approach. Horowitz, Klemelä and Mammen [13] showed that in an additive model with a known identity link function, each additive component can be estimated with the same pointwise normal asymptotic distribution that it would have if the other components were known. Estimation and inference for generalized additive models with known link functions that are not necessarily the identity function have been discussed by Hastie and Tibshirani [11], Linton and Härdle [20], Linton [19], Kauermann and Opsomer [18], Härdle, Huet,

Mammen, and Sperlich [10], Yu, Park and Mammen [43] and Horowitz and Mammen [14]. These models are natural generalizations of generalized linear models (Nelder and Wedderburn [29], Wedderburn [41] and McCullagh and Nelder [28]). Generalized additive models have been put in a larger model framework in Mammen and Nielsen [23]. Generalized additive models with unknown link function have been treated in Horowitz [12] and Horowitz and Mammen [15]. The latter paper generalizes Ichimura's [16] approach for semiparametric single-index models. Coppejans [7] considered a class of additive models that is based on Kolmogorov's theorem on representation of functions of several variables by functions of one variable.

In this paper we will discuss the nonparametric regression model

$$
\begin{aligned}
Y = m\Bigg[ \sum_{l_1=1}^{L_1} m_{l_1}\bigg( \sum_{l_2=1}^{L_2} m_{l_1,l_2}\bigg\{ \cdots \sum_{l_{p-1}=1}^{L_{p-1}} m_{l_1,\ldots,l_{p-1}} \\
\bigg[ \sum_{l_p=1}^{L_p} m_{l_1,\ldots,l_p}(X^{l_1,\ldots,l_p}) \bigg] \bigg\}\bigg) \Bigg] + U,
\end{aligned}
\tag{2}
$$

where $m, m_1, \ldots, m_{L_1,\ldots,L_p}$ are unknown functions and $X^{l_1,\ldots,l_p}$ are one-dimensional elements of a covariate vector $X$, which may be identical for two different indices $(l_1,\ldots,l_p)$. This model is a natural generalization of neural networks where all functions are parametrically specified.

The remainder of the paper is organized as follows. The next two sections discuss the generalized additive model (1). Optimal estimation of the regression function $(x_1,\ldots,x_d) \rightsquigarrow F[m_1(x_1) + \cdots + m_d(x_d)]$ is discussed in Section 2. In Section 3 we show that this result implies that the estimates of the functions $F$ and $m_1,\ldots,m_d$ are rate optimal. Section 4 discusses rate optimal estimation in model (2). Section 5 considers regression quantiles in models (1) and (2). Section 6 presents the results of a simulation study that illustrates the finite-sample performance of our method. Section 7 concludes. The proofs of all results are in Section 8.

**2. Optimal estimation in generalized additive models.** In this section we discuss rate optimal estimation for model (1). We suppose that the response variables $Y_i$ $(i = 1,\ldots,n)$ are given by

$$
Y_i = F[m_1(X_i^1) + \cdots + m_d(X_i^d)] + U_i,
\tag{3}
$$

where $X_i^j$ denotes the $j$th component of the covariate vector $X_i = (X_i^1,\ldots,X_i^d)$, and $X_i$ may be fixed in repeated samples or random. If the covariates are fixed, we assume that the unobserved random variables $U_1,\ldots,U_n$ are independently distributed with $E[U_i] = 0$. If the covariates are random, we assume that $U_1,\ldots,U_n$ are conditionally independent and that $E[U_i|X_i] = 0$.

The functions $F$ and $m_1, \ldots, m_d$ are assumed to belong to a specified class $\mathcal{M}$. $\mathcal{M}$ can be the class of all functions or it can incorporate shape restrictions, such as monotonicity, on some components of $(F, m_1, \ldots, m_d)$.

We estimate $F$ and $m_1, \ldots, m_d$ by penalized least squares. The estimator $(\widehat{F}, \widehat{m}_1, \ldots, \widehat{m}_d)$ minimizes

$$(4) \qquad n^{-1} \sum_{i=1}^{n} \{Y_i - F[m_1(X^1) + \cdots + m_d(X_i^d)]\}^2 + \lambda_n^2 J(F, m_1, \ldots, m_d)$$

over $(F, m_1, \ldots, m_d) \in \mathcal{M}$. Here $J(F, m_1, \ldots, m_d)$ is a penalty term that measures smoothness of order $k$ with $k$ the number of times the functions $F, m_1, \ldots, m_d$ are differentiable. The choice of $J$ is somewhat delicate because we want $J$ to have the same value for all choices of $(F, m_1, \ldots, m_d)$ that result in the same function $(x_1, \ldots, x_d) \to F[m_1(x_1) + \cdots + m_d(x_d)]$. As we discuss below, this can be achieved by the following choice of $J$:

$$J(F, m_1, \ldots, m_d) = J_1^{\nu_1}(F, m_1, \ldots, m_d) + J_2^{\nu_2}(F, m_1, \ldots, m_d),$$

$$J_1(F, m_1, \ldots, m_d) = T_k(F) \left\{ \sum_{j=1}^{d} [T_1^2(m_j) + T_k^2(m_j)] \right\}^{(2k-1)/4},$$

$$J_2(F, m_1, \ldots, m_d) = T_1(F) \left\{ \sum_{j=1}^{d} [T_1^2(m_j) + T_k^2(m_j)] \right\}^{1/4},$$

constants $\nu_1, \nu_2 > 0$ that satisfy $\nu_2 \geq \nu_1$, and

$$T_l^2(f) = \int f^{(l)}(x)^2 \, dx$$

for $0 \leq l \leq k$ and any integrable function $f$. The (possibly random) sequence $(\lambda_n : n = 1, 2, \ldots)$ satisfies conditions that are given in assumption (A5) below. We conjecture that the performance of the estimator does not strongly depend on the choices of the constants $\nu_1, \nu_2$, but we allow here this additional flexibility because a certain choice may simplify the numeric calculation of the estimator.

In fact, the theory that follows does not require $(\widehat{F}, \widehat{m}_1, \ldots, \widehat{m}_d)$ to really minimize (4). It suffices for (4) to differ from its minimum by a term whose size is at most of order $O_P(n^{-2k/(2k+1)})$. In what follows, we will assume that the estimate is chosen so that this holds. This also simplifies the numerical implementation of the estimator. We return to this point below. We call the resulting estimates *approximate* minimizers of (4).

Further normalizing assumptions are needed to identify the functions $(F, m_1, \ldots, m_d)$ in (3). To see this, let $\alpha > 0$ and $\beta = (\beta_1, \ldots, \beta_d) \in \mathbb{R}^d$ be constants. Define

$$(5) \qquad F_{\alpha,\beta}(x) = F[\alpha(x + \beta_1 + \cdots + \beta_d)]$$

and

$$(6) \qquad m_{j,\alpha,\beta}(x) = \alpha^{-1} m_j(x) - \beta_j,$$

for $j = 1, \ldots, d$. Then

$$(7) \qquad F_{\alpha,\beta}[m_{1,\alpha,\beta}(x_1) + \cdots + m_{d,\alpha,\beta}(x_d)] = F[m_1(x_1) + \cdots + m_d(x_d)].$$

Thus, the regression function (conditional mean function of $Y$) is the same for all choices of $\alpha > 0$ and $\beta \in \mathbb{R}^d$. In fact, for a given regression function $H(x) = F[m_1(x_1) + \cdots + m_d(x_d)]$ and under mild regularity conditions, the functions $F$ and $m_1, \ldots, m_d$ are identified up to transformations that correspond to different choices of $\alpha > 0$ and $\beta \in \mathbb{R}^d$. The penalty functionals $J_1$ and $J_2$ are chosen such that they do not depend on the special choice of $\alpha$ and $\beta$. That is,

$$(8) \qquad J_1(F_{\alpha,\beta}, m_{1,\alpha,\beta}, \ldots, m_{d,\alpha,\beta}) = J_1(F, m_1, \ldots, m_d)$$

and

$$(9) \qquad J_2(F_{\alpha,\beta}, m_{1,\alpha,\beta}, \ldots, m_{d,\alpha,\beta}) = J_2(F, m_1, \ldots, m_d)$$

for all $\alpha > 0$ and $\beta \in \mathbb{R}^d$. Therefore, the penalty functionals depend only on the regression function $H(x)$. We will assume that $\mathcal{M}$ is closed under the transformations (5) and (6). See assumption (A3). Then without loss of generality we can assume that $\sum_{j=1}^{d}[T_1^2(m_j) + T_k^2(m_j)] = 1$, and the penalized least squares estimator $(\widehat{F}, \widehat{m}_1, \ldots, \widehat{m}_d)$ can be defined as the minimizer of

$$(10) \qquad \begin{aligned} &\frac{1}{n} \sum_{i=1}^{n} \{Y_i - F[m_1(X_i^1) + \cdots + m_d(X_i^d)]\}^2 \\ &\quad + \lambda_n^2 \left\{ \left[ \int F^{(k)}(z)^2 \, dz \right]^{\nu_1/2} + \left[ \int F'(z)^2 \, dz \right]^{\nu_2/2} \right\} \end{aligned}$$

over all $(F, m_1, \ldots, m_d) \in \mathcal{M}$ with

$$\sum_{j=1}^{d} \left[ \int m_j^{(k)}(x_j)^2 \, dx_j + \int m_j'(x_j)^2 \, dx_j \right] = 1.$$

This norming simplifies the notation when we move to general neural network models in Section 4. But also other scalings are possible and we will use another normalization when we discuss estimation of the additive components and of the link function in Section 3; see (A9) below.

The penalty functionals $J_1$ and $J_2$ contain the $L_2$ norms of the first and $k$th derivatives of $F$ and $m_1, \ldots, m_d$. It can be seen easily that a penalty containing only the $k$th derivatives of these functions will not work here. Consider the extreme case in which $F$ is a linear function. Then $T_k(F) = 0$ and $T_k(m_j)$ can be made arbitrarily small by using the transformations (5) and (6). On the other hand, if $m_1, \ldots, m_d$ are linear functions, then $T_k(m_j) = 0$ for $1 \leq j \leq d$ and $T_k(F)$ can be made arbitrarily small by using the transformations. Therefore, a penalty that depends only on $T_k(F)$ and $T_k(m_1), \ldots, T_k(m_d)$ cannot work because it puts zero penalty on the semiparametric specification in which $F$ or the $m_j$'s are linear.

Our first result states that the regression function $H(x) = F[m_1(x_1) + \cdots + m_d(x_d)]$ can be estimated with rate $n^{-k/(2k+1)}$. This rate is optimal for model (3) with a known link function and unknown additive components under the assumption that the additive components are $k$ times differentiable. Clearly, model (3) is more general, because the link function is unknown. Therefore, this rate is also optimal for (3), and our approach provides a rate optimal estimator.

The rate optimality result needs the following assumptions.

(A1) The covariates $X_i^1, \ldots, X_i^d$ may be fixed in repeated samples or random and take values in a compact subset of $\mathbb{R}$ that, without loss of generality, we take to be $[0, 1]$. The random variables $U_1, \ldots, U_n$ are independent if the covariates are fixed. If the covariates are random, then $U_1, \ldots, U_n$ are conditionally independent given $X_1, \ldots, X_n$.

(A2) The functions $F$ and $m_1, \ldots, m_d$ have $k$ derivatives. Moreover,

$$\int F^{(k)}(x)^2 \, dx < \infty, \qquad \int_0^1 m_j^{(k)}(x)^2 \, dx < +\infty$$

for $j = 1, \ldots, d$. Furthermore, $(F, m_1, \ldots, m_d) \in \mathcal{M}$.

(A3) For all $\alpha > 0$ and $\beta \in \mathbb{R}^d$, if $(G, \mu_1, \ldots, \mu_d) \in \mathcal{M}$, then $(G_{\alpha,\beta}, \mu_{1,\alpha,\beta}, \ldots, \mu_{d,\alpha,\beta}) \in \mathcal{M}$. [For a definition of $G_{\alpha,\beta}, \mu_{1,\alpha,\beta}, \ldots, \mu_{d,\alpha,\beta}$, see (5) and (6).]

(A4) The (conditional) distribution of $U_i$ $(i = 1, \ldots, n)$ has subexponential tails. That is, there are constants $t_U, c_U > 0$ such that

$$\sup_{1 \leq i \leq n} E[\exp(t|U_i|)|X_1, \ldots, X_n] < c_U$$

almost surely for $|t| \leq t_U$. Moreover, $E[U_i | X_1, \ldots, X_n] = 0$ for each $i = 1, \ldots, n$ if the covariates are random, and $E[U_i] = 0$ for each $i = 1, \ldots, n$ if the covariates are fixed in repeated samples.

(A5) $\lambda_n^{-1} = O_p(n^{k/(2k+1)})$ and $\lambda_n = O_p(n^{-k/(2k+1)})$.

These conditions are standard and very weak. In (A1) we assume that the covariates have a compact support to avoid the need of smoothing estimates in the tails of the distribution of $X$. Moreover, a poor rate of convergence for an estimator of one component in the tails could affect the estimator of another component in the center of the distribution of $X$. The (conditional) independence of the $U_i$'s can be weakened to permit martingale difference or mixing sequences of dependent variables. This would complicate the technical analysis and produce a less transparent treatment. Assumption (A2) can be generalized to permit a model that increases with increasing sample size. Again, this would make the theory less transparent and would lead to an estimation procedure in which the sieve model and penalty factors $\lambda_n$ have to be chosen data-adaptively. Assumption (A3) entails no less generality, because $\mathcal{M}$ can always be enlarged to make (A3) hold. Assumption (A4) enables us to use the exponential inequalities needed in empirical process theory. Assumption (A5) allows the possibility that $\lambda_n$ is random. This includes the important case of a data-adaptive choice of $\lambda_n$.

We are now ready to state our first result on rate optimality of our estimator.

THEOREM 2.1. *Let* (A1)–(A5) *hold with* $k \geq 2$. *Then*

$$
\begin{aligned}
(11) \quad & n^{-1} \sum_{i=1}^{n} \{\widehat{F}[\widehat{m}_1(X_i^1) + \cdots + \widehat{m}_d(X_i^d)] \\
& \quad - F[m_1(X_i^1) + \cdots + m_d(X_i^d)]\}^2 = O_p(n^{-2k/(2k+1)})
\end{aligned}
$$

*and*

$$
(12) \qquad J(\widehat{F}, \widehat{m}_1, \ldots, \widehat{m}_d) = O_p(1).
$$

We now state a corollary of Theorem 2.1 for random covariates that satisfy:

(A6) The covariates $X_1, \ldots, X_n$ are independently and identically distributed with distribution $P$.

THEOREM 2.2. *Let* (A1)–(A6) *hold with* $k \geq 2$. *Then*

$$
\begin{aligned}
(13) \quad & \int \{\widehat{F}[\widehat{m}_1(x_1) + \cdots + \widehat{m}_d(x_d)] - F[m_1(x_1) + \cdots + m_d(x_d)]\}^2 P(dx) \\
& \quad = O_p(n^{-2k/(2k+1)})
\end{aligned}
$$

*and* $J(\widehat{F}, \widehat{m}_1, \ldots, \widehat{m}_d) = O_p(1)$.

Up to this point, we have assumed that the penalty factor $\lambda_n$ is the same for all components of $(F, m_1, \ldots, m_d)$. This has been done to simplify the notation. In practice, we can choose a different penalty factor for each component function. To do this, we introduce random factors $\rho_{n,0}, \ldots, \rho_{n,d}$ and modify the penalty functionals $J_1$ and $J_2$ to

$$
J_1(F, m_1, \ldots, m_d) = \rho_{n,0} T_k(F) \left\{ \sum_{j=1}^{d} [T_1^2(m_j) + \rho_{n,j}^2 T_k^2(m_j)] \right\}^{(2k-1)/4}
$$

and

$$J_2(F, m_1, \ldots, m_d) = T_1(F) \left\{ \sum_{j=1}^d [T_1^2(m_j) + \rho_{n,j}^2 T_k^2(m_j)] \right\}^{1/4}.$$

Then Theorems 2.1 and 2.2 hold if $\rho_{n,0}, \ldots, \rho_{n,d} = O_P(1)$ and $\rho_{n,0}^{-1}, \ldots, \rho_{n,d}^{-1} = O_P(1)$.

In this paper we only consider $L_2$ losses. The discussion for sup-norm losses is quite different. Optimal rates differ by different powers of $n$ and not only by a log-term. This can be seen by the construction in the first part of the proof of Theorem 1 in Juditsky, Lepski and Tsybakov [17], which implies that for $d = 2$ and $F$ with $\gamma$ bounded derivatives and $m_1, m_2$ with $\beta$ bounded derivatives up to a logarithmic factor, the order of the optimal rate for sup-norm losses is not faster than $n^{-\gamma/(2\gamma+1+1/\beta)}$. For $\beta = \gamma = 2$, this rate is slower than $n^{-2/5}$. Only if one assumes one more degree of smoothness for $F$ ($\gamma = 3$) does the rate coincide with the optimal $L_2$ rate for $\beta = \gamma = 2$. The basic idea of the construction in Juditsky, Lepski and Tsybakov [17] is to consider testing problems with functions $F$ and $m_2$ both depending on $n$ with shrinking support around zero but with fixed $m_1(x_1) = x_1$. Then for estimating $m_1$ and $m_2$ for $x_1 = x_2 = 0$, only observations $(X_i^1, X_i^2)$ from a local neighborhood around $(0, 0)$ can be used. In Horowitz and Mammen [15] we study pointwise asymptotics of a kernel smoother in an additive model with unknown link under smoothness assumptions $\beta = 2$, $\gamma = 3$ and we show that the pointwise rate $n^{-2/5}$ is achieved.

**3. Optimal estimation of the additive components and link function of a generalized additive model.** Section 2 discussed how well our penalized least squares procedure estimates the conditional mean function, $H(x)$. We now discuss the asymptotic performance of the estimators of the additive components and link function. We make the following additional assumptions.

(A7) The covariates $(X^1, \ldots, X^d)$ have a probability density function $f$ that is bounded away from 0 and $\infty$.

(A8) $F'(z)$ is bounded away from 0 for $z \in \{m_1(x_1) + \cdots + m_d(x_d) : 0 \le x_1, \ldots, x_d \le 1\}$. The additive components $m_j$ are nonconstant for at least two values of $j$ ($1 \le j \le d$).

(A9) The functions $m_1, \ldots, m_d$ and $F$ and their estimates $\widehat{m}_1, \ldots, \widehat{m}_d$ and $\widehat{F}$ are chosen such that

$$\int m_j(x_j) \, dx_j = 0, \qquad \int \widehat{m}_j(x_j) \, dx_j = 0$$

for $j = 1, \ldots, d$ and

$$\sum_{j=1}^d \int m_j(x_j)^2 \, dx_j = 1, \qquad \sum_{j=1}^d \int \widehat{m}_j(x_j)^2 \, dx_j = 1.$$

These are mild conditions. Condition (A7) implies that the $L_2$ norms with respect to the density $f$ and Lebesgue measure are equivalent. This technical point is used in the proof of Theorem 3.2. The assumption that the link function is monotonic is used for identification. All common choices of link functions have this property. The assumption that two additive components are nonconstant is needed for identification. If there were only one nonconstant additive component, say, $m_1$, then it would follow trivially that $F(m_1 + const.)$ does not identify $F$ and $m_1$. Condition (A9) can be always achieved because of (A3) and (A8): Condition (A8) excludes the case that all functions $m_1, \ldots, m_d$ are constant and because of (A3) all functions in $\mathcal{M}$ can be transformed by (5) and (6), at least if not all additive components are constant. Conditions (A8) and (A9) identify the functions $m_1, \ldots, m_d$ and $F$. This can be seen by a simple argument. We state this in the following proposition.

PROPOSITION 3.1. *For continuously differentiable functions* $F \colon \mathbb{R} \to \mathbb{R}$, $m_1 \colon A_1 \to \mathbb{R}, \ldots, m_d \colon A_d \to \mathbb{R}$ *and* $G \colon \mathbb{R} \to \mathbb{R}$, $\mu_1 \colon A_1 \to \mathbb{R}, \ldots, \mu_d \colon A_d \to \mathbb{R}$ *with intervals* $A_1, \ldots, A_d \subset \mathbb{R}$, *we assume that the functions* $m_j$ *are nonconstant for at least two values of* $j$ $(1 \le j \le d)$, $F'(z) > 0$ *for* $z \in \mathbb{R}$,

$$F[m_1(x_1) + \cdots + m_d(x_d)] = G[\mu_1(x_1) + \cdots + \mu_d(x_d)]$$

*for* $x_j \in A_j$, $1 \le j \le d$,

$$\int_{A_j} m_j(x_j)\, dx_j = 0, \qquad \int_{A_j} \mu_j(x_j)\, dx_j = 0$$

*for* $1 \le j \le d$, *and*

$$\sum_{j=1}^{d} \int_{A_j} m_j(x_j)^2\, dx_j = \sum_{j=1}^{d} \int_{A_j} \mu_j(x_j)^2\, dx_j = 1.$$

*Then*

$$m_j(x_j) = \mu_j(x_j)$$

*for* $x_j \in A_j$, $1 \le j \le d$, *and*

$$F(z) = G(z)$$

*for* $z \in \{m_1(x_1) + \cdots + m_d(x_d) \colon x_1 \in A_1, \ldots, x_d \in A_d\}$.

We now state rate-optimality of our estimates of $m_1, \ldots, m_d$ and $F$.

THEOREM 3.2. *Let* (A1)–(A9) *hold with* $k \ge 2$. *Then*

$$(14) \qquad \int_0^1 [\hat{m}_j(x_j) - m_j(x_j)]^2\, dx_j = O_P(n^{-2k/(2k+1)})$$

*and*

$$(15) \quad \int \left\{ \widehat{F}\left[\sum_{j=1}^{d} m_j(x_j)\right] - F\left[\sum_{j=1}^{d} m_j(x_j)\right] \right\}^2 dx = O_P(n^{-2k/(2k+1)}).$$

We now briefly discuss numerical computation of the estimates. We will do this for two approaches. The first is based on B-splines, the second one on smoothing splines. Our estimates are not fully specified because we require only that the penalized least squares objective function be approximately minimized. This leaves some freedom to choose estimates that are best suited to computation. The approach based on B-splines will be used in the simulations below. In this approach we minimize (4) over B-splines $m_1, \ldots, m_d$ and $F$. If the B-splines are of order $k$ and if they use $O(n^{1/(2k+1)})$ knot points, then functions $m_1, \ldots, m_d$ and $F$ that satisfy $T_k(m_1) = O(1), \ldots, T_k(m_d) = O(1)$ and $T_k(F) = O(1)$ can be approximated with an $L_2$ error that is of order $O(n^{-k/(2k+1)})$. This implies that the derivative of $F$ is in supnorm approximated with order o(1) and, thus, $F[m_1(x_1) + \cdots + m_d(x_d)]$ is approximated with order $O(n^{-k/(2k+1)})$. Thus, the minimizer of (4) over B-splines $m_1, \ldots, m_d$ and $F$ is an *approximate minimizer* of (4), as defined in the discussion after (4). The B-spline estimator can be calculated by a backfitting algorithm that alternates between two steps. In one step, $\widehat{F}$ is held fixed at its current value, and a quadratic approximation to the objective function considered as a function of the Fourier coefficients of $m$ is optimized. In the second step, $\widehat{m}$ is held fixed at the value found in the first step, and a new value of $\widehat{F}$ is obtained by optimizing the objective function over the Fourier coefficients of $F$. The first step is an equality-constrained quadratic programming problem that can be solved by the method of Lagrangian multipliers. The second step is an unconstrained quadratic programming problem that can be solved analytically.

The second approach is based on smoothing splines. We will discuss this under the additional assumption that the class $\mathcal{M}$ does not restrict $F$ or one additive component. Condition (A10) makes an assumption for a $j_0$ with $1 \leq j_0 \leq d$.

(A10) For each $(G, \mu_1, \ldots, \mu_d) \in \mathcal{M}$, $(G, \mu_1, \ldots, \mu_{j_0-1}, \mu_{j_0}^*, \mu_{j_0+1}, \ldots, \mu_d) \in \mathcal{M}$ for any function $\mu_{j_0}^* : [0, 1] \rightsquigarrow \mathbb{R}$.

(A11) For each $(G, \mu_1, \ldots, \mu_d) \in \mathcal{M}$, $(G^*, \mu_1, \ldots, \mu_d) \in \mathcal{M}$ for any function $G^* : \mathbb{R} \rightsquigarrow \mathbb{R}$.

THEOREM 3.3. *Let* (A1)–(A8) *hold with* $k \geq 2$.

(i) *Let* (A10) *hold for a* $j_0$ *with* $1 \leq j_0 \leq d$. *Suppose* $(\widetilde{F}, \widetilde{m}_1, \ldots, \widetilde{m}_d)$ *is an approximate minimizer of* (4). *Let* $\bar{m}_{j_0}$ *be chosen among natural splines* $m_{j_0}$ *of order $2k$ with knots $X_1^{j_0}, \ldots, X_n^{j_0}$ so that $(\widetilde{F}, \widetilde{m}_1, \ldots, \widetilde{m}_{j_0-1}, \bar{m}_{j_0}, \widetilde{m}_{j_0+1}, \ldots, \widetilde{m}_d)$ minimizes (4) among $(\widetilde{F}, \widetilde{m}_1, \ldots, \widetilde{m}_{j_0-1}, m_{j_0}, \widetilde{m}_{j_0+1}, \ldots, \widetilde{m}_d)$. Then, $(\widetilde{F}, \widetilde{m}_1, \ldots, \widetilde{m}_{j_0-1}, \bar{m}_{j_0}, \widetilde{m}_{j_0+1}, \ldots, \widetilde{m}_d)$ is also an approximate minimizer of (4) and, therefore has the properties stated in Theorems 2.1, 2.2 and 3.2.*

*(ii) Let* (A11) *hold. Suppose $(\widetilde{F}, \widetilde{m}_1, \ldots, \widetilde{m}_d)$ is an approximate minimizer of (4). Let $\bar{F}$ be chosen among natural splines $F$ of order $2k$ with knots $\widetilde{m}_1(X_1^1) + \cdots + \widetilde{m}_d(X_1^d), \ldots, \widetilde{m}_1(X_n^1) + \cdots + \widetilde{m}_d(X_n^d)$ so that $(\bar{F}, \widetilde{m}_1, \ldots, \widetilde{m}_d)$ minimizes (4) among $(F, \widetilde{m}_1, \ldots, \widetilde{m}_d)$. Then, $(\bar{F}, \widetilde{m}_1, \ldots, \widetilde{m}_d)$ is also an approximate minimizer of (4) and, therefore, has the properties stated in Theorems 2.1, 2.2 and 3.2.*

Natural splines of order $2k$ with knots at the design points arise as minimizers of a penalized least squares criterion for the classical nonparametric regression problem with a one-dimensional regression function and are also called smoothing splines. See, for example, Eubank [8].

We now discuss application of Theorem 3.3 for the case that $\mathcal{M}$ contains all functions. Then (A11) holds and (A10) holds for all $1 \leq j_0 \leq d$. Therefore, repeated application of Theorem 3.3 implies that all estimates, $\widehat{F}$ and $\widehat{m}_1, \ldots, \widehat{m}_d$, can be chosen as natural splines. The computation of the estimates could be done by application of a backfitting algorithm. In each step of the algorithm one estimate $(\widehat{F}, \widehat{m}_1, \ldots,$ or $\widehat{m}_d$, resp.) would be updated. This could be done by using standard smoothing spline software. In the update of $\widehat{m}_1, \ldots, \widehat{m}_d$ the minimization could be approximately solved by linearization.

## 4. Estimation of nonparametric neural network regression.
In this section we discuss rate optimal estimation of the nonparametric neural network model (2). We assume that the response variables $Y_i$ are given by

$$(16) \qquad Y_i = m\left[\sum_{l_1=1}^{L_1} m_{l_1}\left\{\sum_{l_2=1}^{L_2} m_{l_1, l_2}[\cdots m_{l_1, \ldots, l_p}(X_i^{l_1, \ldots, l_p})]\right\}\right] + U_i,$$

where the covariate vector $X_i = (X_i^{l_1, \ldots, l_p} : 1 \leq l_j \leq L_j, 1 \leq j \leq p)$ may be fixed in repeated samples or random. If the covariates are fixed, we assume that the unobserved random variables $U_1, \ldots, U_n$ are independently distributed with $E[U_i] = 0$. If the covariates are random, we assume that the random variables $U_1, \ldots, U_n$ are conditionally independent and that $E[U_i | X_i] = 0$. The functions $(m, m_1, \ldots, m_{L_1, \ldots, L_p})$ are assumed to be contained in a specified class $\mathcal{M}$.

We estimate $(m, m_1, \ldots, m_{L_1, \ldots, L_p})$ by penalized least squares. The penalized least squares estimator $\widehat{m}, \widehat{m}_1, \ldots, \widehat{m}_{L_1, \ldots, L_p}$ minimizes

$$\frac{1}{n}\sum_{i=1}^{n}\left\{Y_i - m\left[\sum_{l_1=1}^{L_1} m_{l_1}\left\{\sum_{l_2=1}^{L_2} m_{l_1, l_2}[\cdots m_{l_1, \ldots, l_p}(X_i^{l_1, \ldots, l_p})]\right\}\right]\right\}^2$$

$$(17) \qquad + \lambda_n^2 J(m)$$

over $(m, m_1, \ldots, m_{L_1,\ldots,L_p}) \in \mathcal{M}$ with

$$J(m) = [T_1^2(m) + cT_k^2(m)]^\nu,$$

$$\sum_{l_1=1}^{L_1} T_1^2(m_{l_1}) + c_{l_1} T_k^2(m_{l_1})$$

$$= \cdots = \sum_{l_p=1}^{L_p} T_1^2(m_{L_1,\ldots,L_{p-1},l_p}) + c_{L_1,\ldots,L_{p-1},l_p} T_k^2(m_{L_1,\ldots,L_{p-1},l_p}) = 1$$

and $\nu, c, c_1, \ldots, c_{L_1,\ldots,L_p} > 0$ constants. It suffices that (17) differs from its minimum by a term that is $O_P(n^{-2k/(2k+1)})$. In what follows, we assume that the estimate is chosen so that this holds.

Our first result states that the regression function $m$ can be estimated with rate $n^{-k/(2k+1)}$, which is optimal for model (16).

THEOREM 4.1. *Let* (A1)–(A5) *hold with* $k \geq 2$, $X_i^1, \ldots, X_i^d$ *replaced by* $X_i^{1,\ldots,1}, \ldots, X_i^{L_1,\ldots,L_p}$ *and* $F, m_1, \ldots, m_d$ *replaced by* $m, m_1, \ldots, m_{L_1,\ldots,L_p}$. *Then*

$$(18) \quad n^{-1} \sum_{i=1}^n \left[ \widehat{m} \left\{ \sum_{l_1=1}^{L_1} \widehat{m}_{l_1}[\cdots \widehat{m}_{l_1,\ldots,l_p}(X_i^{l_1,\ldots,l_p})] \right\} \right.$$
$$\left. - m \left\{ \sum_{l_1=1}^{L_1} m_{l_1}[\cdots m_{l_1,\ldots,l_p}(X_i^{l_1,\ldots,l_p})] \right\} \right]^2 = O_p(n^{-2k/(2k+1)})$$

*and*

$$(19) \qquad J(\widehat{m}) = O_p(1).$$

We now state a corollary of Theorem 4.1 for the case of random covariates.

THEOREM 4.2. *Let* (A1)–(A6) *hold with* $k \geq 2$, *random covariates* $X_i^1, \ldots,$ $X_i^d$ *replaced by* $X_i^{1,\ldots,1}, \ldots, X_i^{L_1,\ldots,L_p}$, *and* $F, m_1, \ldots, m_d$ *replaced by* $m, m_1, \ldots,$ $m_{L_1,\ldots,L_p}$. *Then*

$$(20) \quad \int \left[ \widehat{m} \left\{ \sum_{l_1=1}^{L_1} \widehat{m}_{l_1}[\cdots \widehat{m}_{l_1,\ldots,l_p}(x^{l_1,\ldots,l_p})] \right\} \right.$$
$$\left. - m \left\{ \sum_{l_1=1}^{L_1} m_{l_1}[\cdots m_{l_1,\ldots,l_p}(x^{l_1,\ldots,l_p})] \right\} \right]^2 P(dx) = O_p(n^{-2k/(2k+1)}),$$

*where* $P$ *is the distribution of* $X_i$. *Furthermore,* $J(\widehat{m}) = O_p(1)$.

We conjecture that all functional components can be estimated with the optimal rate $O_p(n^{-k/(2k+1)})$ if (A7) and (A9) hold and $m, \ldots, m_{L_1, \ldots, L_{p-1}}$ have derivatives that are bounded away from 0 and, for all values of $1 \leq l_1 \leq L_1, \ldots, 1 \leq l_{p-1} \leq L_{p-1}$, at least two functions $m_{l_1, \ldots, l_p} : 1 \leq l_p \leq L_p$ are nonconstant. This would be a result that is analogous to Theorem 3.2. Such a result would be less important for neural networks than for generalized additive models. This is because in neural networks one would like to permit two elements of $X$ to be identical, which violates (A7). For example, suppose the regression function is

$$m[m_1[m_{1,1}(x_1) + m_{1,2}(x_2)] + m_1[m_{2,1}(x_1) + m_{2,2}(x_3)]].$$

Arguing as in the proof of Theorem 3.2, one could consistently estimate the partial derivatives $g = m'm_1'm_{1,1}' + m'm_2'm_{2,1}'$, $g_2 = m'm_1'm_{1,2}'$ and $g_3 = m'm_2'm_{2,2}'$. By backfitting, one could fit two functions $h_2(x_1, x_2)$ and $h_3(x_1, x_3)$ such that $g(x_1, x_2, x_3) \approx g_2(x_1, x_2)h_2(x_1, x_2) + g_3(x_1, x_3)h_3(x_1, x_3)$. This would result in estimates of $m_{1,1}'/m_{1,2}'$ and $m_{2,1}'/m_{2,2}'$. Solving, again by backfitting, $\log h_2(x_1, x_2) = h_{2,1}(x_1) + h_{2,2}(x_2)$ and $\log h_3(x_1, x_3) = h_{3,1}(x_1) + h_{3,3}(x_3)$ would give consistent estimates of $m_{1,1}'$, $m_{2,1}'$, $m_{1,2}'$ and $m_{2,2}'$. It is clear that it is very hard to establish the conditions under which this approach would result in a consistent estimate. It would be even harder to show that this approach can be used to get rate optimal estimates of the functions $m$, $m_1$, $m_2$, $m_{1,1}$, $m_{1,2}$, $m_{2,1}$ and $m_{2,2}$.

**5. Regression quantiles.** The estimation approach of this paper can be extended to $M$-functionals other than least squares. In this section we will discuss quantile estimation. We consider again model (1) or (16), but now we choose $0 < \alpha < 1$ and we assume that the (conditional) $\alpha$-quantile of $U_i$ is equal to 0 (and not the conditional mean). We define $u_\alpha(z) = \alpha z - zI[z \leq 0]$, where $I$ is the indicator function. Define penalized regression quantiles as the functions that minimize [up to a term of order $O_P(n^{-2k/(2k+1)})$]

$$(21) \quad \frac{1}{n}\sum_{i=1}^{n} u_\alpha\{Y_i - \widehat{F}[\widehat{m}_1(X_i^1) + \cdots + \widehat{m}_d(X_i^d)]\} + \lambda_n^2 J(\widehat{F}, \widehat{m}_1, \ldots, \widehat{m}_d)$$

or

$$(22) \quad \frac{1}{n}\sum_{i=1}^{n} u_\alpha\left[Y_i - \widehat{m}\left\{\sum_{l_1=1}^{L_1} \widehat{m}_{l_1}[\cdots \widehat{m}_{l_1, \ldots, l_p}(X_i^{l_1, \ldots, l_p})]\right\}\right] + \lambda_n^2 J(\widehat{m}).$$

The penalty terms are as defined in Sections 2 and 4. Make the following assumption.

(A4′) The function $E[u_\alpha(U_i - \mu)|X_1, \ldots, X_n]$ almost surely has a unique minimum at $\mu = 0$. Furthermore, for some $\varepsilon > 0$ and all $0 \leq a \leq \varepsilon$, it holds that

$$\inf_{1 \leq i \leq n} P(0 \leq U_i \leq a) \geq \varepsilon a$$

TABLE 1
*Performance of $\widehat{m}_1$, $\widehat{m}_2$ and $\widehat{F}$ for different values of $n$ and $\lambda$*

| **n** | **$\lambda$** | **$\widehat{m}_1$** | **$\widehat{m}_2$** | **$\widehat{F}$** |
|---|---|---|---|---|
| 400 | 0.05 | 0.030 | 0.029 | 0.0040 |
|     | 0.10 | 0.029 | 0.024 | 0.0039 |
|     | 0.15 | 0.026 | 0.029 | 0.0048 |
| 900 | 0.05 | 0.023 | 0.018 | 0.0030 |
|     | 0.10 | 0.017 | 0.015 | 0.0027 |
|     | 0.15 | 0.025 | 0.017 | 0.0036 |

and

$$\inf_{1 \le i \le n} P(-a \le U_i \le 0) \ge \varepsilon a$$

almost surely.

THEOREM 5.1. *Let the conditions of Theorem 2.1, Theorem 2.2, Theorem 3.2, Theorem 4.1 or Theorem 4.2 hold with $(\mathrm{A}4')$ in place of $(\mathrm{A}4)$. Then the conclusions of the corresponding theorem hold for the estimators defined in (21) or (22).*

**6. Simulation results.** We carried out a small simulation study with $Y = F[m_1(X^1) + m_2(X^2)] + U$, where $F$ is the identity function, $m_1(x) = \sin(\pi x)$, $m_2(x) = \Phi(3x)$, $\Phi$ is the standard normal distribution, and $U \sim \mathrm{N}(0,1)$. The values of $(X^1, X^2)$ are the grid $(i/(n^{1/2} + 1), j/(n^{1/2} + 1))$, $i, j = 1, \ldots, n^{1/2}$, where $n$ is the sample size. The penalty term $J$ is defined with $\nu_1 = \nu_2 = 1$. We used the B-spline approach described in Section 3. The estimates of $m_1$, $m_2$ and $F$ are B-splines with four knots. There are 500 Monte Carlo replications in each simulation.

Table 1 shows the empirical integrated mean-square errors of $\widehat{m}_1$, $\widehat{m}_2$ and $\widehat{F}$ at three different values of the penalty parameter, $\lambda$.

The simulation results with $\lambda = 0.10$ are shown graphically in Figures 1 and 2. The wiggles in the estimates of $\widehat{m}_2$ are due to variance, not bias. The 4-knot spline fits the true $m_2$ very well. In the simulations our estimators show a very reliable performance.

**7. Conclusions and extensions.** In this paper we have proposed an estimation approach for a general class of nested regression models. The basic idea is to use the following property of compositions of functions belonging to certain smoothness classes: if the same entropy rate applies for all smoothness classes, then the same entropy rate also applies to the class of the composition of the functions. In our setting, the function classes are subsets of additive Sobolev classes. The results could be extended easily to other smoothness classes as long as entropy rates with respect to the supremum norm are available. Examples are additive Sobolev classes of functions with higher-dimensional arguments. Another point that needs exploration is the case in which smoothness classes with different entropy rates enter into the model. It would be interesting to check whether each component's convergence rate is the one corresponding to the entropy rate of its smoothness class. In particular, for parametric components it would be important to check whether the component can be estimated with rate $n^{-1/2}$. Furthermore, we conjecture that the resulting estimate is efficient. Such a result has been proved in Mammen and van de Geer [27] for a partial linear model with a known link function. There, penalized quasi-likelihood estimation is used for the nonparametric components. Another extension would be to apply our results for other classes of $M$ estimators.

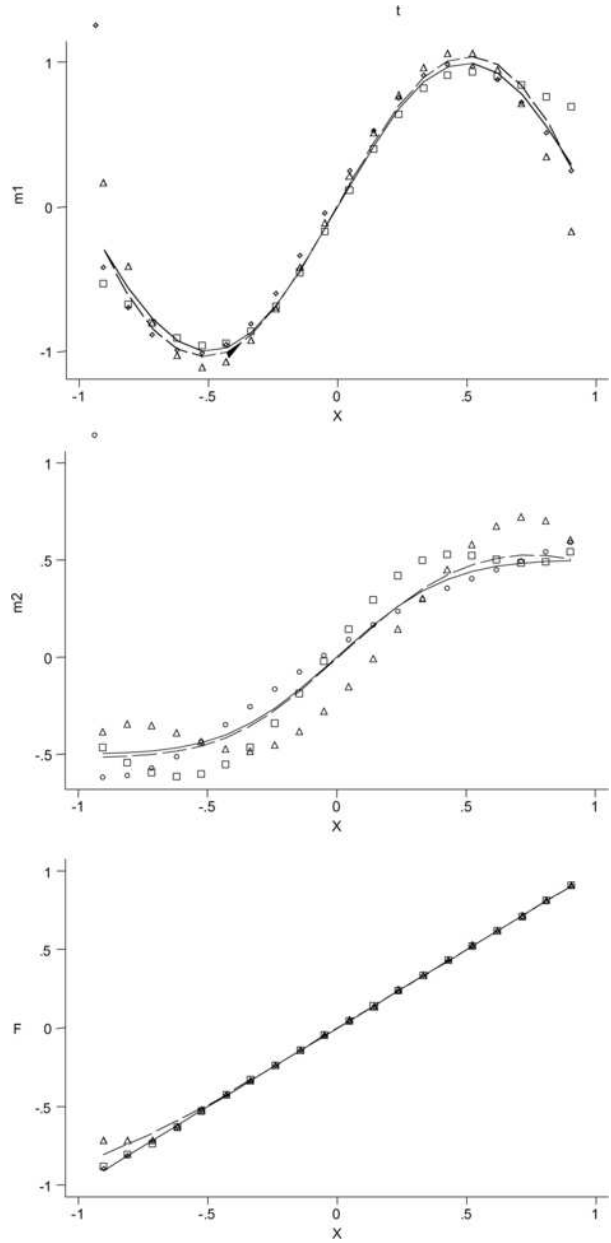

Fig. 1. *Performance of $\widehat{m}_1$ (upper plot), $\widehat{m}_2$ (middle plot) and $\widehat{F}$ (lower plot) with $n = 400$. The solid line is the true function; the dashed line is average of 500 estimates; circles, squares and triangles, respectively, denote the estimates at the 25th, 50th and 75th percentiles of the IMSE.*

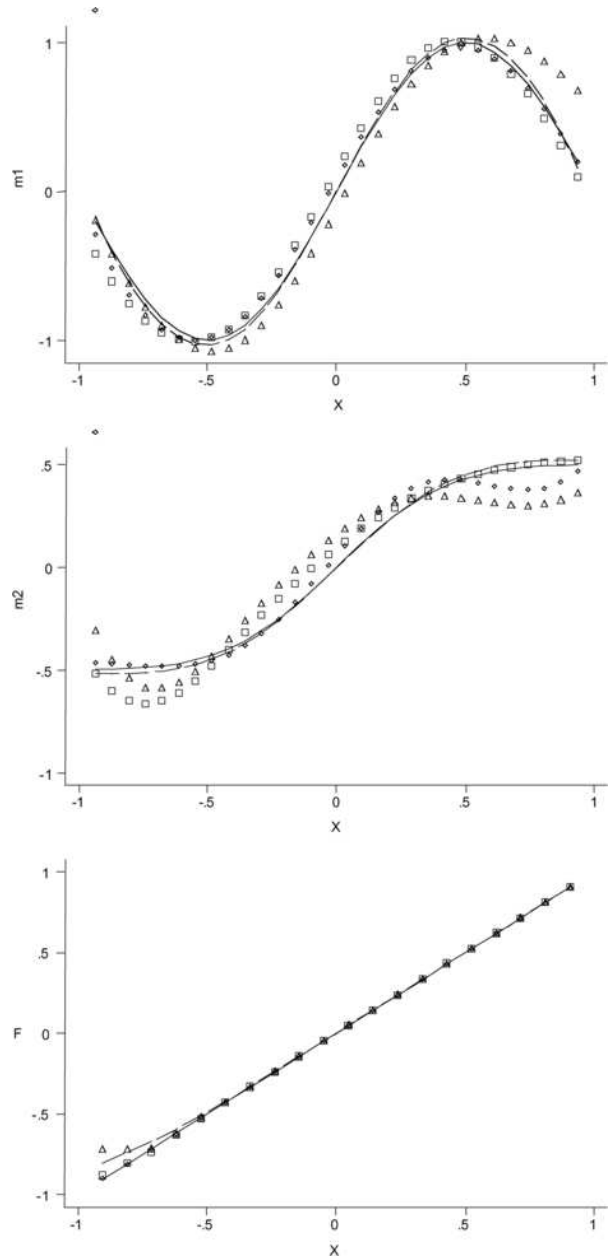

Fig. 2.   *As Figure 1 but with* $n = 900$.

## 8. Proofs.

PROOF OF THEOREM 2.1. For a constant $c > 0$ consider the class of functions

$$\mathcal{G} = \left\{ F[m_1(x_1) + \cdots + m_d(x_d)] \colon |F(z)| \le c \text{ for } |z| \le d, m_j(0) = 0 \right.$$

$$\text{for } j = 1, \ldots, d, \ \sum_{j=1}^{d} \int_0^1 m_j^{(k)}(x)^2 \, dx + \sum_{j=1}^{d} \int_0^1 m_j'(x)^2 \, dx = 1,$$

$$\left. J(F, m_1, \ldots, m_d) \le 1 \right\}.$$

First we will argue that, for a constant $C_K$,

$$(23) \qquad H_B(\delta, \mathcal{G}, \| \cdot \|_\infty) \le C_k \delta^{-1/k}$$

for $\delta > 0$. Here, $\| \cdot \|_\infty$ denotes the supremum norm. Furthermore, $H_B(\delta, \mathcal{G}, \| \cdot \|_\infty)$ denotes the $\delta$-entropy with bracketing for the class $\mathcal{G}$ w.r.t. the sup norm $\| \cdot \|_\infty$. This means that $\exp(H_B)$ is the smallest number $N$ for which there exist pairs of functions $(g_1^L, g_1^U), \ldots, (g_N^L, g_N^U)$ in $\mathcal{G}$ with the following property. For each $g \in \mathcal{G}$ there exists $1 \le j \le N$ with $g_j^L \le g \le g_j^U$ and $\| g_j^U - g_j^L \|_\infty \le \delta$. Such a set of tuples is also called a $\delta$-cover with bracketing.

This entropy bound follows from the following classical entropy bound on Sobolev classes (see Birman and Solomjak [3] and van de Geer [40]):

$$(24) \qquad H_B\left(\delta, \left\{ g \colon [0,1] \to \mathbb{R} \colon \|g\|_\infty \le 1, \right. \right.$$

$$\left. \left. \int g^{(k)}(x)^2 \, dx \le 1 \right\}, \| \cdot \|_\infty \right) \le C \delta^{-1/k}$$

for a constant $C > 0$. We now show how (23) follows from (24). From (24) one gets for the class of additive functions

$$\mathcal{G}_{\text{add}} = \left\{ m_1(x_1) + \cdots + m_d(x_d) : \sum_{j=1}^d \int_0^1 m_j^{(k)}(x)^2 \, dx + \sum_{j=1}^d \int_0^1 m_j'(x)^2 \, dx \le 1, \right.$$
$$\left. m_j(0) = 0 \right\}$$

with a constant $C' > 0$

$$H_B(\delta, \mathcal{G}_{\text{add}}, \| \cdot \|_\infty) \le C' \delta^{-1/k}. \tag{25}$$

We use here that $\int_0^1 m_j'(x)^2 \, dx \le 1$ and $m_j(0) = 0$ implies that $\|m_j\|_\infty \le 1$.

Consider now a function $F[m_1(x_1) + \cdots + m_d(x_d)]$ that is an element of $\mathcal{G}$. Suppose that $m_1, \ldots, m_d$ are chosen such that $m_j(0) = 0$ for $j = 1, \ldots, d$ and $\sum_{j=1}^d \int_0^1 m_j^{(k)}(x)^2 \, dx + \sum_{j=1}^d \int_0^1 m_j'(x)^2 \, dx = 1$. For such a representation $J(F, m_1, \ldots, m_d) \le 1$ implies $\int F^{(k)}(z)^2 \, dz \le 1$. Because $|F(z)| \le c$ for $|z| \le d$, this implies $|F'(z)| \le C'$ for $|z| \le d$ with a constant $C'$. This can be seen, for example, by application of the interpolation inequality; see (42). Consider now a $\delta$-cover with bracketing $(g_1^L, g_1^U), \ldots, (g_N^L, g_N^U)$ of $\mathcal{G}_{\text{add}}$. Consider a fixed function $F$ with $0 \le |F'| \le C'$. Then $[F(g_1^L) - C'\delta, F(g_1^U) + C'\delta], \ldots, [F(g_L^L) - C'\delta, F(g_L^U) + C'\delta]$ is a $(2C'\delta)$-cover with bracketing of $F(\mathcal{G}_{\text{add}})$. By a slight extension of this argument, we get (23).

We now apply Theorem 10.2 in van de Geer [40] with the modifications discussed before the theorem. This theorem implies (11) and (12). We now verify the assumptions of Theorem 10.2 in van de Geer [40]. We have to check for $\varepsilon > 0$ that, with probability larger than $1 - \varepsilon$, the function $\widehat{H}(x) = \widehat{F}^*[\widehat{m}_1(x_1) + \cdots + \widehat{m}_d(x_d)]$ is an element of $\mathcal{G}$ if $c$ is chosen large enough. Here the function $\widehat{F}^*$ is defined as $\widehat{F}^*(z) = \widehat{F}(z) \, / (1 + J(\widehat{F}, \widehat{m}_1, \ldots, \widehat{m}_d))$. W.l.o.g. we can assume that

$$\sum_{j=1}^d \int \widehat{m}_j^{(k)}(x)^2 \, dx + \int \widehat{m}_j'(x)^2 \, dx = 1, \tag{26}$$

$$\widehat{m}_j(0) = 0 \qquad \text{for } 1 \le j \le d. \tag{27}$$

It can be easily checked that $J(\widehat{F}^*, \widehat{m}_1, \ldots, \widehat{m}_d) \le 1$. Thus, for the proof of $\widehat{H} \in \mathcal{G}$, it remains to check that

$$\sup_{|z| \le d} |\widehat{F}^*(z)| = O_P(1). \tag{28}$$

We now show (28). Equations (26) and (27) imply that

$$\sup_{0 \le x_1, \ldots, x_d \le 1} |\widehat{m}_1(x_1) + \cdots + \widehat{m}_d(x_d)| \le d.$$

Furthermore, because of $J(\widehat{F}^*, \widehat{m}_1, \ldots, \widehat{m}_d) \le 1$, these equations imply that $\int \widehat{F}^{*\prime}(z)^2\,dz \le 1$. This shows that

$$\sup_{|z|,|z'| \le d} |\widehat{F}^*(z') - \widehat{F}^*(z)| \le 2d.$$

We now show that

(29) $$\inf_{|z| \le d} |\widehat{F}^*(z)| = O_P(1).$$

The last two bounds imply (28). Thus, it remains to show (29). For the proof of (29), note first that by definition of $\widehat{F}, \widehat{m}_1, \ldots, \widehat{m}_d$ the following inequality holds with $\widetilde{F}(z) \equiv \overline{Y} = n^{-1} \sum_{i=1}^n Y_i$ and $Z_i = \widehat{m}_1(X_i^1) + \cdots + \widehat{m}_d(X_i^d)$:

$$\frac{1}{n} \sum_{i=1}^n \{Y_i - \widehat{F}[Z_i]\}^2 \le \frac{1}{n} \sum_{i=1}^n \{Y_i - \widehat{F}[Z_i]\}^2 + \lambda_n^2 J(\widehat{F}, \widehat{m}_1, \ldots, \widehat{m}_d)$$

$$\le \frac{1}{n} \sum_{i=1}^n \{Y_i - \widetilde{F}[Z_i]\}^2 + \lambda_n^2 J(\widetilde{F}, \widehat{m}_1, \ldots, \widehat{m}_d)$$

$$= \frac{1}{n} \sum_{i=1}^n \{Y_i - \overline{Y}\}^2$$

$$= O_P(1).$$

This implies $\inf_{|z| \le d} |\widehat{F}(z)| = O_P(1)$ because of

$$\left| \overline{Y} - \frac{1}{n} \sum_{i=1}^n \widehat{F}(Z_i) \right|^2 = \left| \frac{1}{n} \sum_{i=1}^n Y_i - \widehat{F}(Z_i) \right|^2 \le \frac{1}{n} \sum_{i=1}^n \{Y_i - \widehat{F}[Z_i]\}^2.$$

Claim (29) now follows because of $|\widehat{F}^*| \le |\widehat{F}|$. $\quad\square$

PROOF OF THEOREM 2.2. For the proof of Theorem 2.2, it remains to show (13). This claim immediately follows from Lemma 5.16 in van de Geer [40]. $\quad\square$

PROOF OF PROPOSITION 3.1. Without loss of generality, we assume that the functions $m_1$ and $m_2$ are nonconstant. Then there exist $x_1^* \in A_1$ and $x_2^* \in A_2$ with $m_1'(x_1^*) \ne 0$ and $m_2'(x_2^*) \ne 0$. For $H(x) = F[m_1(x_1) + \cdots + m_d(x_d)] = G[\mu_1(x_1) + \cdots + \mu_d(x_d)]$, we get that $\frac{\partial}{\partial x_1} H(x) \ne 0$ if $x_1 = x_1^*$ and $\frac{\partial}{\partial x_2} H(x) \ne 0$ if $x_2 = x_2^*$. For $x_1 \in A_1, \ldots, x_d \in A_d$, put $x^* = (x_1^*, x_2, \ldots, x_d)'$ and $x^{**} = (x_1, x_2^*, x_3, \ldots, x_d)'$. Then for $2 \le j \le d$,

$$\frac{m_j'(x_j)}{m_1'(x_1^*)} = \frac{\frac{\partial}{\partial x_j} H(x^*)}{\frac{\partial}{\partial x_1} H(x^*)} = \frac{\mu_j'(x_j)}{\mu_1'(x_1^*)}.$$

Because of $\int_{A_j} m_j(x_j)\,dx_j = \int_{A_j} \mu_j(x_j)\,dx_j = 0$, this gives for $2 \leq j \leq d$

$$(30) \qquad \frac{m_j(x_j)}{m_1'(x_1^*)} = \frac{\mu_j(x_j)}{\mu_1'(x_1^*)}.$$

Using partial derivatives of $H$ at $x = x^{**}$, we get

$$(31) \qquad \frac{m_1(x_1)}{m_2'(x_2^*)} = \frac{\mu_1(x_1)}{\mu_2'(x_2^*)}.$$

Equation (31) implies that $m_1'(x_1^*)/m_2'(x_2^*) = \mu_1'(x_1^*)/\mu_2'(x_2^*)$. This shows that (30) holds for $1 \leq j \leq d$. Because of $\sum_{j=1}^{d} \int_{A_j} m_j(x_j)^2\,dx_j = \sum_{j=1}^{d} \int_{A_j} \mu_j(x_j)^2\,dx_j$, this implies the statements of the proposition. $\square$

PROOF OF THEOREM 3.2. We first show (14). Put $\widehat{H}(x_1,\ldots,x_d) = \widehat{F}[\widehat{m}_1(x_1) + \cdots + \widehat{m}_d(x_d)]$ and $H(x_1,\ldots,x_d) = F[m_1(x_1) + \cdots + m_d(x_d)]$. We write $H_j = \partial_{x_j} H$, $\widehat{H}_j = \partial_{x_j} \widehat{H}$, $H_{i,j} = \partial_{x_j} \partial_{x_i} H$ and $\widehat{H}_{i,j} = \partial_{x_j} \partial_{x_i} \widehat{H}$ for the partial derivatives of $H$ and $\widehat{H}$.

For $1 \leq j \leq d$, define $\widetilde{m}_j(x_j) = \widehat{\gamma}^{-1}[\widehat{m}_j(x_j) - \widehat{m}_j(0)]$ with $\widehat{\gamma}^2 = \sum_{j=1}^{d} \int \widehat{m}_j^{(k)}(x)^2\,dx + \int \widehat{m}_j'(x)^2\,dx$. Furthermore, choose $\widetilde{F}$ so that $\widetilde{F}[\widetilde{m}_1(x_1) + \cdots + \widetilde{m}_d(x_d)] = \widehat{F}[\widehat{m}_1(x_1) + \cdots + \widehat{m}_d(x_d)]$.

Then $\widetilde{m}_1,\ldots,\widetilde{m}_d$ satisfy (26) and (27) with $\widehat{m}_j$ replaced by $\widetilde{m}_j$ and we have that

$$(32) \qquad \int_0^1 \widetilde{m}_j(x_j)^2\,dx_j = O_P(1),$$

$$(33) \qquad \int_0^1 \widetilde{m}_j^{(k)}(x_j)^2\,dx_j = O_P(1)$$

for $j = 1,\ldots,d$. Note also that $\widetilde{m}_j(0) = 0$ by definition.

By Sobolev embedding results (see, e.g., Section VI.7 in Yosida [42] or Oden and Reddy [32]), the bounds (32) and (33) give

$$(34) \qquad \sup_{x_j \in [0,1]} |\widetilde{m}_j^{(l)}(x_j)| = O_P(1)$$

for $j = 1,\ldots,d$ and $0 \leq l \leq k-1$. We now derive a similar bound for the link function $F$.

From (29) and (12) one gets that $\inf_{|z| \leq d} |\widetilde{F}(z)| = O_P(1)$. From Theorem 2.1 we get that

$$(35) \qquad \int \widetilde{F}^{(k)}(z)^2\,dz = O_P(1).$$

By application of the Sobolev embedding, this shows that

$$(36) \qquad \sup_{|z| \leq d} |\widetilde{F}^{(l)}(z)| = O_P(1)$$

for $0 \le l \le k-1$.

The rest of the proof is divided into several steps.

STEP 1. In this step we argue that

$$(37) \qquad \int_{[0,1]^d} \widehat{H}^*(x)^2 \, dx = O_P(1),$$

where $\widehat{H}^*$ is a partial derivative of $\widehat{H}$ of order $k$. The integral $\int \widehat{H}^*(x)^2 \, dx$ can be easily bounded by a sum of integrals over products of derivatives of $\widetilde{F}$, $\widetilde{m}_1$ or $\ldots$ or $\widetilde{m}_d$, respectively. Most summands can be easily bounded by using (33)–(36). One summand needs a little bit more care, namely,

$$\int_{[0,1]^d} \widetilde{F}^{(k)} (\widetilde{m}_1(x_1) + \cdots + \widetilde{m}_d(x_d))^2 \widetilde{m}'_{i_1}(x_{i_1})^2 \cdot \cdots \cdot \widetilde{m}'_{i_k}(x_{i_k})^2 \, dx.$$

This term arises when $H^*$ is a partial derivative w.r.t. $x_{i_1}, \ldots, x_{i_k}$. Up to a factor that is stochastically bounded, this integral is equal to

$$(38) \qquad \int_{[0,1]^d} \widetilde{F}^{(k)} (\widetilde{m}_1(x_1) + \cdots + \widetilde{m}_d(x_d))^2 \widetilde{m}'_{i_1}(x_{i_1})^2 \, dx;$$

see also (34). We now apply that for two functions $g : [0,1] \to [a,b]$, $f : [a,b] \to \mathbb{R}$ with $a < b$ the following inequality holds:

$$(39) \qquad \int_0^1 f[g(y)]^2 g'(y)^2 \, dy \le 2 \left[ \int_0^1 g''(y)^2 \, dy \right]^{1/2} \int_a^b f(z)^2 \, dz.$$

By using (39) with $f = \widetilde{F}^{(k)}$ and $g = \widetilde{m}_1 + const.$, one can easily check that the integral in (38) is bounded by

$$O_P(1) \cdot \int_{-d}^d \widetilde{F}^{(k)}(z)^2 \, dz.$$

This quantity is stochastically bounded because of (35). For the proof of (37), it remains to prove (39). For the proof of this inequality, we denote for $u < v$ by $k(u,v)$ the number of crossings of the interval $[u,v]$ by the function $g'$. It can be easily checked that

$$\int_{I_{u,v}} |g''(y)| \, dy \ge (v-u) k(u,v),$$

where

$$I_{u,v} = \{ y \in [0,1] : u < g'(y) \le v \}.$$

Choose now $c_i = 2^{-i}$. The claim (39) now follows from

$$\int_0^1 f[g(y)]^2 g'(y)^2 \, dy = \int_{\{y : g'(y) \ne 0\}} f[g(y)]^2 g'(y)^2 \, dy$$

$$= \sum_{i=-\infty}^{\infty} \int_{I_{c_i,c_{i-1}}} f[g(y)]^2 g'(y)^2 \, dy$$

$$+ \sum_{i=-\infty}^{\infty} \int_{I_{-c_{i-1},-c_i}} f[g(y)]^2 g'(y)^2 \, dy$$

$$\leq \sum_{i=-\infty}^{\infty} [k(c_i, c_{i-1}) + k(-c_{i-1}, -c_i)] c_{i-1} \int_a^b f(z)^2 \, dz$$

$$\leq \sum_{i=-\infty}^{\infty} \frac{c_{i-1}}{c_{i-1} - c_i} \int_{I_{c_i,c_{i-1}} \cup I_{-c_{i-1},-c_i}} |g''(y)| \, dy \int_a^b f(z)^2 \, dz$$

$$\leq 2 \int_0^1 |g''(y)| \, dy \int_a^b f(z)^2 \, dz$$

$$\leq 2 \Big[ \int_0^1 g''(y)^2 \, dy \Big]^{1/2} \int_a^b f(z)^2 \, dz.$$

STEP 2. We now show that

$$\int [\widehat{H}_j(x_1, \ldots, x_d) - H_j(x_1, \ldots, x_d)]^2 \, dx = O_P(n^{-2(k-1)/(2k+1)}), \tag{40}$$

$$\int [\widehat{H}_{i,j}(x_1, \ldots, x_d) - H_{i,j}(x_1, \ldots, x_d)]^2 \, dx = O_P(n^{-2(k-2)/(2k+1)}) \tag{41}$$

for $1 \leq i, j \leq d$.

For the proof of these claims, we make use of the interpolation inequality of Agmon [2]; see also van de Geer ([40], Lemma 10.8) and Mammen and Thomas-Agnan [26]. This inequality states that for a function $g : \mathbb{R} \to \mathbb{R}$ and a real number $\theta > 0$ it holds that

$$\int g^{(l)}(x)^2 \, dx \leq c\theta^{-2l} \int g(x)^2 \, dx + c\theta^{2(k-l)} \int g^{(k)}(x)^2 \, dx \tag{42}$$

for a constant $c$ and $1 \leq l \leq k$. The claims (40) and (41) follow from the bound on $\widehat{H} - H$ in Theorem 2.1, (37) and the interpolation inequality.

STEP 3. According to (A7), two additive functions are not constant. W.l.o.g. we assume that this is the case for the first two functions. Then there exist constants $0 \leq a_1 < b_1 \leq 1$ and $0 \leq a_2 < b_2 \leq 1$ with

$$\inf_{a_j \leq x_j \leq b_j} |m_j'(x_j)| > 0 \qquad \text{for } j = 1, 2.$$

In this step we show that uniformly for $0 \leq x_1 \leq 1$ it holds that

$$\widehat{\rho}\widetilde{m}_1'(x_1) = \rho m_1'(x_1) + o_P(1), \tag{43}$$

where

$$\widehat{\rho} = \int_{a_1}^{b_1} \frac{1}{\widetilde{m}_1'(x_1)} \, dx_1,$$

$$\rho = \int_{a_1}^{b_1} \frac{1}{m_1'(x_1)} \, dx_1.$$

For the proof of (43) note first that (40)–(41) imply that there exist random $0 \le x_3^*, \ldots, x_d^* \le 1$, $a_2 \le x_2^* \le b_2$ with

$$(44) \quad \begin{aligned} &\int [\widehat{H}_j(x_1, x_2^*, \ldots, x_d^*) - H_j(x_1, x_2^*, \ldots, x_d^*)]^2 \, dx_1 \\ &\quad = O_P(n^{-2(k-1)/(2k+1)}), \end{aligned}$$

$$(45) \quad \begin{aligned} &\int [\widehat{H}_{j,1}(x_1, x_2^*, \ldots, x_d^*) - H_{j,1}(x_1, x_2^*, \ldots, x_d^*)]^2 \, dx_1 \\ &\quad = O_P(n^{-2(k-2)/(2k+1)}) \end{aligned}$$

for $j = 1$ and $j = 2$. We now argue that for a (random) function $\Delta : [0,1] \to \mathbb{R}$ the following implication holds. If $\int \Delta'(u)^2 \, du = O_P(1)$ and $\int \Delta(u)^2 \, du = o_P(1)$, then it holds that $\sup |\Delta(u)| = o_P(1)$. This implication can be easily verified by using that $\int \Delta'(u)^2 du = O_P(1)$ implies that

$$\sup_{0 \le u,v \le 1} \frac{|\Delta(u) - \Delta(v)|}{|u - v|^{1/2}} = O_P(1).$$

The latter implication follows by application of an embedding theorem (see Adams [1], page 97) or directly by a simple calculation.

We now apply this result for $j = 1$ and $j = 2$ with $\Delta(x_1) = \widehat{H}_j(x_1, x_2^*, \ldots, x_d^*) - H_j(x_1, x_2^*, \ldots, x_d^*)$. This gives

$$\sup_{0 \le x_1 \le 1} |\widehat{H}_j(x_1, x_2^*, \ldots, x_d^*) - H_j(x_1, x_2^*, \ldots, x_d^*)| = o_P(1).$$

We now apply this expansion and make use of the fact that $|m_1'|(u)$ for $(u \in [a_1, b_1])$, $|m_2'|(u)$ for $(u \in [a_2, b_2])$ and $F'$ are bounded away from zero and from infinity. We get the following expansions that hold uniformly for $0 \le x_1 \le 1$ and $a_1 \le x_1' \le b_1$:

$$\begin{aligned} \frac{\widetilde{m}_2'(x_2^*)}{\widetilde{m}_1'(x_1')} &= \frac{\widehat{H}_2(x_1', x_2^*, \ldots, x_d^*)}{\widehat{H}_1(x_1', x_2^*, \ldots, x_d^*)} \\ &= \frac{H_2(x_1', x_2^*, \ldots, x_d^*)}{H_1(x_1', x_2^*, \ldots, x_d^*)} + o_P(1) \\ &= \frac{m_2'(x_2^*)}{m_1'(x_1')} + o_P(1), \end{aligned}$$

$$\frac{\widetilde{m}_1'(x_1)}{\widetilde{m}_2'(x_2^*)} = \frac{m_1'(x_1)}{m_2'(x_2^*)} + o_P(1).$$

This implies that uniformly for $0 \le x_1 \le 1$ and $a_1 \le x_1' \le b_1$,

$$\frac{\widetilde{m}_1'(x_1)}{\widetilde{m}_1'(x_1')} = \frac{m_1'(x_1)}{m_1'(x_1')} + o_P(1).$$

Claim (43) now follows by integrating both sides of the last equality w.r.t. $x_1'$.

STEP 4.   In this step we show that for $2 \le j \le d$ and for random sequences $\delta_{j,n}$,

$$(46) \qquad \int_0^1 |\widetilde{m}_j(x_j) - \widehat{\rho}^{-1}\rho[m_j(x_j) - m_j(0)] - \delta_{j,n}|^2 \, dx_j = O_P(n^{-2k/(2k+1)}).$$

For the proof we note first that (40) and the bound on $\widehat{H} - H$ in Theorem 2.1 imply that there exist random numbers $0 \le x_2^*, \ldots, x_{j-1}^*, x_{j+1}^*, \ldots, x_d^* \le 1$ with

$$(47) \qquad \begin{aligned} &\int [\widehat{H} - H]^2(x_1, x_2^*, \ldots, x_{j-1}^*, x_j, x_{j+1}^*, \ldots, x_d^*) \, dx_1 \, dx_j \\ &\qquad = O_P(n^{-2k/(2k+1)}), \end{aligned}$$

$$(48) \qquad \begin{aligned} &\int [\widehat{H}_1 - H_1]^2(x_1, x_2^*, \ldots, x_{j-1}^*, x_j, x_{j+1}^*, \ldots, x_d^*) \, dx_1 \, dx_j \\ &\qquad = O_P(n^{-2(k-1)/(2k+1)}), \end{aligned}$$

$$(49) \qquad \begin{aligned} &\int [\widehat{H}_j - H_j]^2(x_1, x_2^*, \ldots, x_{j-1}^*, x_j, x_{j+1}^*, \ldots, x_d^*) \, dx_1 \, dx_j \\ &\qquad = O_P(n^{-2(k-1)/(2k+1)}). \end{aligned}$$

In the following calculations of this step we fix the random vector $(x_2^*, \ldots, x_{j-1}^*, x_{j+1}^*, \ldots, x_d^*)$ and, for simplicity of notation, we write $f(x_1, x_j)$ instead of $f(x_1, x_2^*, \ldots, x_{j-1}^*, x_j, x_{j+1}^*, \ldots, x_d^*)$ for the functions $f = H, H_1, H_j, \widehat{H}, \widehat{H}_1$ or $\widehat{H}_j$, respectively. We now use that

$$(50) \qquad \begin{aligned} \int_0^{x_j} du_j \int_{a_1}^{b_1} dx_1 \frac{\widehat{H}_j(x_1, u_j)}{\widehat{H}_1(x_1, u_j)} &= \int_{a_1}^{b_1} \frac{1}{\widetilde{m}_1'(x_1)} \, dx_1 [\widetilde{m}_j(x_j) - \widetilde{m}_j(0)] \\ &= \widehat{\rho}\widetilde{m}_j(x_j) \end{aligned}$$

and that

$$\int_0^{x_j} du_j \int_{a_1}^{b_1} dx_1 \frac{H_j(x_1, u_j)}{H_1(x_1, u_j)}$$

(51)
$$= \int_{a_1}^{b_1} \frac{1}{m_1'(x_1)} \, dx_1 [m_j(x_j) - m_j(0)]$$
$$= \rho[m_j(x_j) - m_j(0)].$$

Furthermore, we make use of the expansion

$$\frac{1}{\widehat{H}_1} = \frac{1}{H_1} - \frac{\widehat{H}_1 - H_1}{H_1^2} + (\widehat{H}_1 - H_1)^2 \frac{1}{H_1^2 \widehat{H}_1}.$$

This gives the expansion

$$\widetilde{m}_j(x_j) - \widehat{\rho}^{-1} \rho[m_j(x_j) - m_j(0)]$$

(52)
$$= \widehat{\rho}^{-1} \int_0^{x_j} du_j \int_{a_1}^{b_1} dx_1 \left[ \frac{\widehat{H}_j(x_1, u_j)}{\widehat{H}_1(x_1, u_j)} - \frac{H_j(x_1, u_j)}{H_1(x_1, u_j)} \right]$$

$$= \widehat{\rho}^{-1} \int_0^{x_j} du_j \int_{a_1}^{b_1} dx_1 \left[ \frac{\widehat{H}_j - H_j}{H_1} - \frac{(\widehat{H}_1 - H_1)(\widehat{H}_j - H_j)}{H_1^2} \right.$$

$$\left. - H_j \frac{\widehat{H}_1 - H_1}{H_1^2} + \widehat{H}_j (\widehat{H}_1 - H_1)^2 \frac{1}{H_1^2 \widehat{H}_1} \right] (x_1, u_j).$$

Because of (34), it holds that $\widehat{\rho}^{-1} = O_P(1)$. This bound together with (48)–(49) implies for the second term in (52) the bound

(53)
$$\widehat{\rho}^{-1} \int_0^{x_j} du_j \int_{a_1}^{b_1} dx_1 \left[ \frac{(\widehat{H}_1 - H_1)(\widehat{H}_j - H_j)}{H_1^2} \right] (x_1, u_j)$$

$$= O_P(n^{-(2k-2)/(2k+1)})$$

$$= O_P(n^{-k/(2k+1)}).$$

For estimating the last term in (52), we use that

$$\sup_{a_1 \leq x_1 \leq b_1, 0 \leq x_j \leq 1} \left| \widehat{\rho}^{-1} \frac{\widehat{H}_j}{\widehat{H}_1}(x_1, x_j) \right| = \sup_{a_1 \leq x_1 \leq b_1, 0 \leq x_j \leq 1} \left| \widehat{\rho}^{-1} \frac{\widetilde{m}_j'(x_j)}{\widetilde{m}_1'(x_1)} \right| = O_P(1),$$

because of (40) and because $\inf_{a_1 \leq x_1 \leq b_1} |\widehat{\rho} \widetilde{m}_1'(x_1)| > c$ with probability tending to one for a constant $c > 0$ small enough. The latter fact follows directly from (43) and (A6). With this bound, we get for the last term in (52)

(54)
$$\widehat{\rho}^{-1} \int_0^{x_j} du_j \int_{a_1}^{b_1} dx_1 \left[ \widehat{H}_j (\widehat{H}_1 - H_1)^2 \frac{1}{H_1^2 \widehat{H}_1} \right] (x_1, u_j)$$

$$= O_P(1) \int_0^{x_j} du_j \int_{a_1}^{b_1} dx_1 (\widehat{H}_1 - H_1)^2 (x_1, u_j)$$

$$= O_P(n^{-k/(2k+1)}),$$

where again (48) was applied. Using (52)–(54), we get uniformly for $0 \leq x_j \leq 1$ that

$$
\begin{aligned}
\widetilde{m}_j(x_j) &- \widehat{\rho}^{-1} \rho[m_j(x_j) - m_j(0)] \\
&= \widehat{\rho}^{-1} \int_0^{x_j} du_j \int_{a_1}^{b_1} dx_1 \left[ \frac{\widehat{H}_j - H_j}{H_1} \right](x_1, u_j) \\
&\quad - \widehat{\rho}^{-1} \int_0^{x_j} du_j \int_{a_1}^{b_1} dx_1 H_j \left[ \frac{\widehat{H}_1 - H_1}{H_1^2} \right](x_1, u_j) + O_P(n^{-k/(2k+1)}) \\
&= T_1(x_j) + T_2(x_j) + O_P(n^{-k/(2k+1)}).
\end{aligned}
$$

We now apply partial integration for the first term. This gives

$$
\begin{aligned}
T_1(x_j) &= \widehat{\rho}^{-1} \int_{a_1}^{b_1} dx_1 \left[ \frac{1}{H_1(x_1, u_j)} \{ \widehat{H}(x_1, u_j) - H(x_1, u_j) \} \right]_{u_j=0}^{u_j=x_j} \\
&\quad + \widehat{\rho}^{-1} \int_{a_1}^{b_1} dx_1 \int_0^{x_j} du_j \left[ \frac{H_{1,j}(x_1, u_j)}{H_1(x_1, u_j)^2} \{ \widehat{H}(x_1, u_j) - H(x_1, u_j) \} \right] \\
&= -\widehat{\rho}^{-1} \int_{a_1}^{b_1} dx_1 \left[ \frac{1}{H_1(x_1, 0)} \{ \widehat{H}(x_1, 0) - H(x_1, 0) \} \right] \\
&\quad + \widehat{\rho}^{-1} \int_{a_1}^{b_1} dx_1 \left[ \frac{1}{H_1(x_1, x_j)} \{ \widehat{H}(x_1, x_j) - H(x_1, x_j) \} \right] \\
&\quad + \widehat{\rho}^{-1} \int_{a_1}^{b_1} dx_1 \int_0^{x_j} du_j \left[ \frac{H_{1,2}(x_1, u_j)}{H_1(x_1, u_j)^2} \{ \widehat{H}(x_1, u_j) - H(x_1, u_j) \} \right] \\
&= \delta_{j,1,n} + g_{j,1,n}(x_j) + g_{j,2,n}(x_j),
\end{aligned}
$$

with a real random sequence $\delta_{j,1,n}$. The random functions $g_{j,1,n}$ and $g_{j,2,n}$ satisfy

$$
\int g_{j,1,n}(x_j)^2 \, dx_j = O_P(n^{-2k/(2k+1)}),
$$

$$
\int g_{j,2,n}(x_j)^2 \, dx_j = O_P(n^{-2k/(2k+1)}).
$$

Similarly, we get that

$$
T_2(x_j) = \delta_{j,2,n} + g_{j,3,n}(x_j) + g_{j,4,n}(x_j)
$$

with a real random sequence $\delta_{j,2,n}$ and random functions $g_{j,3,n}$ and $g_{j,4,n}$ that satisfy

$$
\int g_{j,3,n}(x_j)^2 \, dx_j = O_P(n^{-2k/(2k+1)}),
$$

$$
\int g_{j,4,n}(x_j)^2 \, dx_j = O_P(n^{-2k/(2k+1)}).
$$

This shows that, for $2 \leq j \leq d$,

$$\int_0^1 |\widetilde{m}_j(x_j) - \widehat{\rho}^{-1} \rho[m_j(x_j) - m_j(0)] - \delta_{j,1,n} - \delta_{j,2,n}|^2 \, dx_j = O_P(n^{-2k/(2k+1)}).$$

This implies (46).

STEP 5.   In this step we show that there exists a random sequence $\delta_{1,n}$ such that

(55)    $$\int_0^1 |\widetilde{m}_1(x_1) - \widehat{\rho}^{-1} \rho[m_1(x_1) - m_1(0)] - \delta_{1,n}|^2 \, dx_1 = O_P(n^{-2k/(2k+1)}).$$

For this purpose we choose a function $s \colon [0,1] \to \mathbb{R}$ that has a continuous derivative and satisfies $s(0) = s(1) = 0$ and $\int s(x_2) m_2'(x_2) \, dx_2 = 1$. Put $w(x_2) = s(x_2) m_2'(x_2) \widetilde{m}_2'(x_2)$. One can easily check that

$$\widetilde{m}_1(x_1) = \int_0^{x_1} du_1 \int_0^1 dx_2 \, w(x_2) \frac{\widehat{H}_1(u_1, x_2)}{\widehat{H}_2(u_1, x_2)},$$

where $\widehat{H}$ is defined as in the last step for $j = 2$. We define

$$m_1^*(x_1) = \int_0^{x_1} du_1 \int_0^1 dx_2 \, w(x_2) \frac{H_1(u_1, x_2)}{H_2(u_1, x_2)}.$$

Proceeding as above, one can show that there exists a random sequence $\delta_{1,n}$ such that

$$\int_0^1 |\widetilde{m}_1(x_1) - m_1^*(x_1) - \delta_{1,n}|^2 \, dx_1 = O_P(n^{-2k/(2k+1)}).$$

In particular, the proof makes use of the following facts: $\sup_{x_2} |w(x_2)| = O_P(1)$, $\sup_{x_1, \ldots, x_d} |w(x_2)[\widehat{H}_1/\widehat{H}_2](x_1, x_2, \ldots, x_d)| = O_P(1)$ and $\int w'(x_2)^2 \, dx_2 = O_P(1)$. For the proof of (55), it remains to show that

$$\int_0^1 |m_1^*(x_1) - \widehat{\rho}^{-1} \rho[m_1(x_1) - m_1(0)]|^2 \, dx_1 = O_P(n^{-2k/(2k+1)}).$$

Because of

$$m_1^*(x_1) = [m_1(x_1) - m_1(0)] \int_0^1 s(x_2) \widetilde{m}_2'(x_2) \, dx_2,$$

this follows from

$$\int_0^1 s(x_2) \widetilde{m}_2'(x_2) \, dx_2$$

$$= s(x_2) \widetilde{m}_2(x_2)|_0^1 - \int_0^1 s'(x_2) \widetilde{m}_2(x_2) \, dx_2$$

$$= -\int_0^1 s'(x_2) \widehat{\rho}^{-1} \rho m_2(x_2) \, dx_2 + O_P(n^{-k/(2k+1)})$$

$$= \int_0^1 s(x_2)\widehat{\rho}^{-1}\rho m_2'(x_2)\,dx_2 + O_P(n^{-k/(2k+1)})$$

$$= \widehat{\rho}^{-1}\rho + O_P(n^{-k/(2k+1)}).$$

Here, (46) for $j = 2$ was used. Thus, (55) is proved.

STEP 6.   In this step we show that

$$(56) \qquad\qquad\qquad \widehat{\rho} = O_P(1).$$

Arguing as above, we find $x_2^*, \ldots, x_d^*$ such that

$$\sup_{0 \le x_1 \le 1} |H(x_1, x_2^*, \ldots, x_d^*) - \widehat{H}(x_1, x_2^*, \ldots, x_d^*)| = o_P(1).$$

The claim now directly follows from

$$0 < \inf_z F'(z)|m_1(b_1) - m_1(a_1)|$$

$$\le |F[m_1(b_1) + m_2(x_2^*) + \cdots + m_d(x_d^*)]$$

$$\quad - F[m_1(a_1) + m_2(x_2^*) + \cdots + m_d(x_d^*)]|$$

$$= |\widetilde{F}[\widetilde{m}_1(b_1) + \widetilde{m}_2(x_2^*) + \cdots + \widetilde{m}_d(x_d^*)]$$

$$\quad - \widetilde{F}[\widetilde{m}_1(a_1) + \widetilde{m}_2(x_2^*) + \cdots + \widetilde{m}_d(x_d^*)]| + o_P(1)$$

$$\le \sup_z \widetilde{F}'(z)|\widetilde{m}_1(b_1) - \widetilde{m}_1(a_1)| + o_P(1)$$

$$= O_P(1)\widehat{\rho}^{-1}\rho|m_1(b_1) - m_1(a_1)| + o_P(1).$$

STEP 7.   In this step we show claim (14).

Using $\int m_j(x_j)\,dx_j = \int \widehat{m}_j(x_j)\,dx_j = 0$, the definition of $\widetilde{m}_j$, (46) and (55), we get for $1 \le j \le d$

$$(57) \qquad \int_0^1 |\widehat{\gamma}^{-1}\widehat{m}_j(x_j) - \widehat{\rho}^{-1}\rho m_j(x_j)|^2\,dx_j = O_P(n^{-2k/(2k+1)}).$$

Here we have used that for a function $w$ it holds that $\int [w(x) - \int w(u)\,du]^2\,dx \le \int w(x)^2\,dx$.

We now use that for a constant $\alpha > 0$ and for two functions $w_1$ and $w_2$ with $\int w_1(x)^2\,dx = \int w_2(x)^2\,dx = 1$, it holds that $\int [\alpha w_1(x) - \alpha^{-1}w_2(x)]^2\,dx \le \int [w_1(x) - w_2(x)]^2\,dx$. This shows

$$\widehat{\gamma}^{-1}\widehat{\rho}^{-1}\rho \sum_{j=1}^d \int_0^1 |\widehat{m}_j(x_j) - m_j(x_j)|^2\,dx_j$$

$$\le \widehat{\gamma}^{-1}\widehat{\rho}^{-1}\rho \sum_{j=1}^d \int_0^1 \left| \sqrt{\widehat{\gamma}^{-1}\widehat{\rho}\rho^{-1}}\,\widehat{m}_j(x_j) - \frac{1}{\sqrt{\widehat{\gamma}^{-1}\widehat{\rho}\rho^{-1}}}m_j(x_j) \right|^2\,dx_j$$

$$= O_P(n^{-2k/(2k+1)}).$$

Furthermore, because of (56) and (34), we have $\widehat{\rho}\rho^{-1} = O_P(1)$ and $\widehat{\rho}^{-1}\rho = O_P(1)$. With (57), this gives $\widehat{\gamma}^{-2} = \widehat{\rho}^{-2}\rho^2 + o_P(1)$ and, thus, $\widehat{\gamma} = O_P(1)$. Therefore, the last inequality implies (14).

STEP 8. It remains to show (15). From (36) with $l = 1$ and (57), we get

$$(58) \qquad \sup_{|z| \leq d} |\widehat{F}'(z)| = O_P(1).$$

Claim (15) immediately follows from (14), (58) and Theorem 2.2. $\square$

PROOF OF THEOREM 3.3. We give a proof of part (i). Part (ii) follows by similar arguments. Suppose that $(\widetilde{F}, \widetilde{m}_1, \ldots, \widetilde{m}_d)$ is an approximate minimizer of (4) over $\mathcal{M}$ with $\widetilde{m}_{j_0}$ not necessarily a natural spline. Define $\bar{m}_{j_0}$ as minimizer of $\int \mu^{(k)}(u)^2 \, du$ under the constraint $\mu(X_i^{j_0}) = \widetilde{m}_{j_0}(X_i^{j_0})$ for $1 \leq i \leq n$. The function $\bar{m}_{j_0}$ is a natural spline of order $2k$ with knots $X_1^{j_0}, \ldots, X_n^{j_0}$; see, for example, Eubank [8]. We show that

$$(59) \qquad \int \bar{m}'_{j_0}(u)^2 \, du = O_P(1).$$

This immediately shows that $(\widetilde{F}, \widetilde{m}_1, \ldots, \widetilde{m}_{j_0-1}, \bar{m}_{j_0}, \widetilde{m}_{j_0+1}, \ldots, \widetilde{m}_d)$ is an approximate minimizer of (4) over $\mathcal{M}$ and, thus, it implies the statement of Theorem 3.2(i). It remains to show (59). This follows from

$$(60) \qquad \int \widehat{\Delta}'(u)^2 \, du = O_P(1),$$

with $\widehat{\Delta}(u) = \widetilde{m}_{j_0}(u) - \bar{m}_{j_0}(u)$. For the proof of (60) note that, by the Sobolev embedding theorem, one can write $\widehat{\Delta}(z) = \widehat{\Delta}_1(z) + \widehat{\Delta}_2(z)$ with

$$\widehat{\Delta}_1(z) = \sum_{j=1}^{k} \widehat{\beta}_j z^{j-1},$$

and $|\widehat{\Delta}_2(z)| \leq [\int \widehat{\Delta}^{(k)}(z)^2 \, dz]^{1/2} = O_P(1)$; see, for example, Oden and Reddy [32]. Because of $\widehat{\Delta}(X_i^{j_0}) = 0$ for $1 \leq i \leq n$, we get that $\widehat{\beta}_1, \ldots, \widehat{\beta}_k = O_P(1)$. This implies $\int \widehat{\Delta}(u)^2 \, du = O_P(1)$. Now (60) follows from the interpolation inequality (42). $\square$

PROOF OF THEOREMS 4.1 AND 4.2. The theorems follow by similar arguments as in the proofs of Theorems 2.1 and 2.2. $\square$

PROOF OF THEOREM 5.1. The proof of the quantile version of Theorem 2.1 and Theorem 4.1 follows along the same lines as in the old proofs. For the necessary modifications to apply empirical process theory, see van de Geer [39] and Chapter 12 in van de Geer [40]. Note, for example, that (for $\alpha = 1/2$) condition (A4′) restates (12.22) and (12.23) in van de Geer [40]. Compare also Exercise 12.4 in van de Geer [40]. The quantile versions of Theorem 2.2, Theorem 3.2 and Theorem 4.2 directly follow from the new versions of Theorem 2.1 and Theorem 4.1 by the same arguments as in the proofs of their old versions. □

**Acknowledgments.** We are very grateful for the helpful, constructive and very detailed comments of two referees. We thank Sasha Tsybakov for referring us to the lower bound in Juditsky, Lepski and Tsybakov [17].

DEPARTMENT OF ECONOMICS
NORTHWESTERN UNIVERSITY
EVANSTON, ILLINOIS 60208-2600
USA
E-MAIL: joel-horowitz@northwestern.edu

DEPARTMENT OF ECONOMICS
UNIVERSITY OF MANNHEIM
L7, 3–5
68131 MANNHEIM
GERMANY
E-MAIL: emammen@rumms.uni-mannheim.de